\documentclass[11pt,reqno]{amsart}

\usepackage[T1]{fontenc}
\usepackage{lmodern}
\usepackage{microtype}
\usepackage{amsmath,amssymb,mathtools}
\usepackage{graphicx}
\usepackage{tikz}
\usepackage{placeins}
\usepackage{hyperref}
\hypersetup{
  hidelinks,
  pdftitle={A transitive C\textasciicircum(1+\unichar{945}) surface diffeomorphism with positive entropy has at most one measure of maximal entropy},
  pdfauthor={Wanshan Lin and Xueting Tian},
  pdfsubject={Smooth ergodic theory},
  pdfkeywords={Measure of maximal entropy, surface diffeomorphism, topological transitivity, nonuniform hyperbolicity}
}

\numberwithin{equation}{section}
\newtheorem{theorem}{Theorem}[section]
\newtheorem{lemma}[theorem]{Lemma}
\newtheorem{proposition}[theorem]{Proposition}

\DeclareMathOperator{\diam}{diam}
\newcommand{\htop}{h_{\mathrm{top}}}
\newcommand{\doi}[1]{\href{https://doi.org/#1}{\nolinkurl{doi:#1}}}

\title[Maximal entropy measures for surface diffeomorphisms]{A transitive
  \(C^{1+\alpha}\) surface diffeomorphism\\
  with positive entropy has at most one\\
  measure of maximal entropy}
\author{Wanshan Lin}
\address{School of Mathematical Sciences, Fudan University, Shanghai 200433, People's Republic of China}
\email{wanshanlin@fudan.edu.cn}

\author{Xueting Tian}
\address{School of Mathematical Sciences, Fudan University, Shanghai 200433, People's Republic of China}
\email{xuetingtian@fudan.edu.cn}
\date{}
\begin{document}
\begin{abstract}
Let \(0<\alpha\le1\). On a closed surface, every topologically transitive
\(C^{1+\alpha}\) diffeomorphism with positive topological entropy
has at most one measure of maximal entropy.
\end{abstract}
\maketitle

\section{Introduction}
Let \(f\colon M\to M\) be a diffeomorphism of a \emph{closed surface},
that is, a compact smooth two-dimensional manifold without boundary. An invariant
Borel probability measure \(\mu\) is a \emph{measure of maximal entropy} if
\(h_\mu(f)=\htop(f)\). Sarig proved that when \(f\) is \(C^{1+\alpha}\)
and \(\htop(f)>0\), at most countably many ergodic measures attain this
maximum \cite{Sarig}.

Buzzi, Crovisier and Sarig proved that a topologically transitive
\(C^\infty\) surface diffeomorphism of positive topological entropy has a
unique measure of maximal entropy~\cite{BCS}. For \(C^r\) surface
diffeomorphisms with \(1<r<\infty\), Burguet's theorem~\cite{Burguet}
gives existence when \(\htop(f)>\lambda_{\min}(f)/r\)
(see also \cite[Remark~1.2 and Theorem~A]{BLY} for nonintegral \(r\)).
Under this bound, topological transitivity implies uniqueness~\cite{BCS}. Here
\[
 \lambda_{\min}(f)=\min\{\lambda(f),\lambda(f^{-1})\},\qquad
 \lambda(g)=\lim_{n\to\infty}\frac1n\log\sup_{x\in M}\|D_xg^n\|.
\]
Without the strict entropy bound, existence can fail, as Buzzi's \(C^r\)
examples on closed surfaces show~\cite{Buzzi}. Buzzi, Crovisier and Sarig
announce such examples at the threshold
\(\htop(f)=\lambda_{\min}(f)/r>0\)~\cite[Section~1.8(c)]{BCS}.

For \(1\le r<\infty\), Buzzi, Crovisier and Sarig asked whether a
topologically transitive \(C^r\) surface diffeomorphism of positive topological
entropy can have several maximizing measures
\cite[Question~1]{BCS}. For every \(1<r<\infty\), we show that positive
topological entropy suffices to rule out such multiplicity.

\begin{theorem}\label{thm:main}
Fix \(0<\alpha\le1\) and let \(f\) be a \(C^{1+\alpha}\) diffeomorphism
of a closed surface. If \(f\) is topologically transitive and
\(\htop(f)>0\), then \(f\) has at most one measure of maximal entropy.
\end{theorem}

\begin{samepage}
Equivalently, under the same regularity and transitivity assumptions,
any nonatomic ergodic measure of maximal entropy whose Lyapunov exponents
have opposite signs is the unique measure of maximal entropy.
The existence of such a measure implies
\(\htop(f)>0\) by Katok's criterion~\cite[Corollary~4.2]{Katok}.
\par\end{samepage}

\begin{samepage}
The conclusion can fail at zero entropy or in dimension three.
Windsor~\cite[Section~2]{Windsor} constructed a minimal \(C^\infty\)
diffeomorphism of \(\mathbb T^2\) with exactly two ergodic invariant probability
measures, and Banerjee and Kunde~\cite[Theorem~D]{BK} gave a real-analytic example.
Each has zero topological entropy, since positive entropy would force a hyperbolic
periodic point~\cite[Corollary~4.3]{Katok}, contrary to minimality.
In each case both ergodic measures therefore have maximal entropy.
\par\end{samepage}

On \(\mathbb T^3\), there exist robustly mixing diffeomorphisms with positive
topological entropy and at least two ergodic measures of maximal
entropy~\cite[Theorem~2]{RHRHTU}.

\subsection*{Organization of this paper}
Section~\ref{sec:branches} constructs return branches and horseshoes in small
rectangles, and Section~\ref{sec:passages} uses topological transitivity to
connect these rectangles by passages of fixed duration. Section~\ref{sec:count}
completes the proof by comparing two entropy bounds for the labels of
synchronized returns.

\subsection*{AI Declaration}
We used ChatGPT (GPT-6, including Astra) to help develop and review steps in
the proof. It helped us show why one passage time works for every horizontal
crossing of the source rectangle and check the two entropy bounds
for synchronized returns obtained by conditioning on the split and on the
return time, respectively. We also used it to locate references,
check the hypotheses of cited results, revise the English text, and prepare
the diagrams and their captions. We reviewed the AI-assisted material,
revised it as needed, and take full responsibility for the article.

\section{Return branches in a nice domain}\label{sec:branches}

A \emph{return rectangle} is a nice domain \(R=\Gamma_{pq}\) of
Climenhaga, Luzzatto and Pesin, a closed topological disk with periodic
saddles \(p,q\) at opposite corners, \emph{vertical sides} in
\(W^s_{\mathrm{loc}}(p),W^s_{\mathrm{loc}}(q)\), and \emph{horizontal sides}
in \(W^u_{\mathrm{loc}}(p),W^u_{\mathrm{loc}}(q)\).
Its associated integer \(T(R)\) is an even common multiple of the periods
of \(p,q\). At positive multiples of \(T(R)\), the stable sides avoid
\(\operatorname{int}R\) under forward iteration and the unstable sides do so
under backward iteration \cite[Definition~1.16 and Remark~1.17]{CLP}.

\begin{figure}[htbp]
\centering
\begin{tikzpicture}[x=1cm,y=1cm,font=\footnotesize,
    crossing/.style={line width=1.05pt,line cap=round}]
  \fill[black!12] (1.10,0) rectangle (2.30,1.66);
  \draw[black!60] (1.10,0) -- (1.10,1.66) (2.30,0) -- (2.30,1.66);
  \draw (0,0) rectangle (3.47,1.66);
  \node at (3.17,1.40) {$R$};
  \node at (1.69,.23) {$S$};
  \draw (.00,.69) .. controls (.39,.71) and (.45,1.32) .. (.72,1.25)
                    .. controls (.88,1.21) and (.96,.99) .. (1.10,.95);
  \draw[crossing] (1.10,.95) .. controls (1.39,.87) and (1.52,.55) .. (1.78,.59)
                             .. controls (2.00,.63) and (2.12,.93) .. (2.30,.97);
  \draw (2.30,.97) .. controls (2.63,1.05) and (2.95,.66) .. (3.47,.78);
  \fill (1.10,.95) circle (1.15pt) (2.30,.97) circle (1.15pt);
  \node[above] at (.64,1.23) {$\gamma$};
  \node at (1.68,.96) {$\gamma_0$};
  \draw[->] (3.72,.83) -- node[above] {$f^n|_S$} (5.02,.83);
  \begin{scope}[xshift=5.27cm]
    \fill[black!12] (0,.42) rectangle (3.47,1.23);
    \draw[black!60] (0,.42) -- (3.47,.42) (0,1.23) -- (3.47,1.23);
    \draw (0,0) rectangle (3.47,1.66);
    \node at (3.17,1.43) {$R$};
    \draw[crossing] (0,.96) .. controls (.72,.95) and (1.02,.67) .. (1.61,.71)
                            .. controls (2.20,.75) and (2.72,1.08) .. (3.47,1.01);
    \fill (0,.96) circle (1.15pt) (3.47,1.01) circle (1.15pt);
    \node at (1.64,1.02) {$f^n(\gamma_0)$};
    \node at (2.88,.64) {$f^n(S)$};
  \end{scope}
\end{tikzpicture}
\caption{Both panels show \(R\), with stable sides on the left and right
and unstable sides at the top and bottom. A horizontal crossing \(\gamma\)
has a subpath \(\gamma_0\) joining the vertical sides of \(S\).
\(f^n(\gamma_0)\) crosses \(R\) inside the image strip \(f^n(S)\).}
\label{fig:return-strip}
\end{figure}
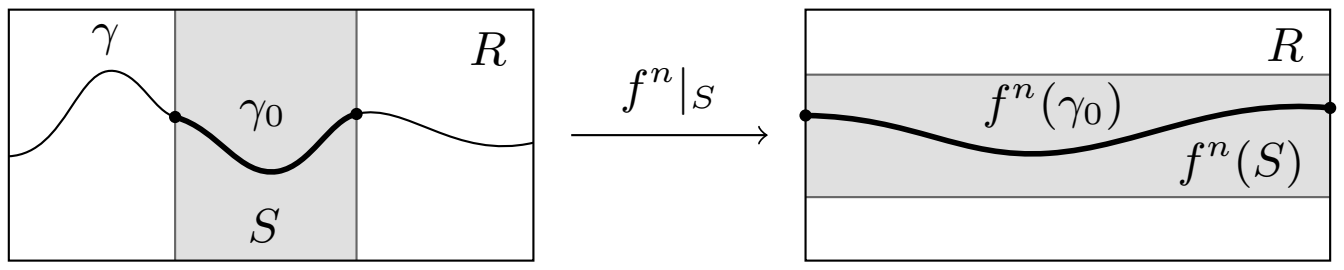

A \emph{horizontal crossing} is a continuous path in \(R\) joining its two vertical
sides. A \emph{vertical strip} is a subrectangle of \(R\) whose top and bottom sides
lie on the corresponding sides of \(R\). A \emph{horizontal strip} is defined
similarly using the left and right sides. A \emph{return branch} of \emph{order} \(n\) in
Lemma~\ref{lem:branches} is a restriction \(f^n|_S\) mapping its vertical
\emph{source strip} \(S\subset R\) onto its horizontal \emph{image strip}
\(f^n(S)\subset R\) (Figure~\ref{fig:return-strip}).
Intermediate visits to \(R\) are allowed.

An ergodic measure is of \emph{saddle type} if its two Lyapunov exponents
have opposite signs.

\begin{lemma}\label{lem:branches}
Let \(\mu\) be a nonatomic ergodic measure of saddle type for a
\(C^{1+\alpha}\) diffeomorphism \(f\) of a closed surface.
For each \(r>0\) there exist a return rectangle \(R\) with \(\diam R<r\),
a Borel set \(A\subset\operatorname{int}R\) with \(\mu(A)>0\), a positive
integer \(T\), and \(B,\kappa>0\) such that
\begin{enumerate}
\item for each \(n\in T\mathbb N\), the set \(A\cap f^{-n}A\) is covered
by at most countably many pairwise disjoint closed vertical strips
\(S\subset R\), each mapped by \(f^n\) onto a horizontal strip in \(R\),
\item every horizontal crossing of \(R\) has a subpath in each such \(S\)
whose image under \(f^n\) crosses \(R\) horizontally, and
\item for any sequence of \(m\ge1\) branches \(f^{n_i}|_{S_i}\) from the
families in (1), with cumulative times \(t_0=0\) and \(t_i=n_1+\cdots+n_i\),
the intersection \(S=\bigcap_{i=1}^m f^{-t_{i-1}}(S_i)\) is a nonempty
closed vertical strip. With \(n=t_m\), the concatenated
branch \(f^n|_S\) maps \(S\) onto a horizontal strip in \(R\) and sends
the two vertical sides of \(S\) into opposite vertical sides of \(R\). Moreover,
\begin{equation}\label{eq:thin}
  \diam f^k(S)\leq B\bigl(e^{-\kappa k}+e^{-\kappa(n-k)}\bigr)
  \qquad (0\leq k\leq n).
\end{equation}
\end{enumerate}
\end{lemma}

\begin{proof}
At the level \(\ell\) supplied by \cite[Theorem~C]{CLP}, choose a nice
regular set \(E\subset R=\Gamma_{pq}\) with \(\mu(E)>0\) and \(\diam R\)
below both \(r\) and the threshold for that level in \cite[Theorem~D]{CLP}.
Each \(z\in E\) has full-length local stable and
unstable arcs \(W_z^s,W_z^u\) joining the horizontal and vertical sides of
\(R\), respectively~\cite[p.~994]{CLP}. These sides lie in stable or unstable
manifolds of periodic saddles, so recurrence and nonatomicity give \(\mu(\partial R)=0\).
Choose a compact \(A\subset E\cap\operatorname{int}R\) with \(\mu(A)>0\)
and put \(T=2T(R)\), so \(Df^T\) preserves orientation on the stable and
unstable lines at \(p,q\).

For \(n\in T\mathbb N\) and \(x\in A\cap f^{-n}A\), put \(y=f^nx\in A\subset E\).
Since \(y\in f^n(W_x^s)\cap W_y^u\), this is an almost return to \(E\)
\cite[Definition~4.13]{CLP}. By \cite[Theorem~D and its proof]{CLP},
\(x\) lies in a closed vertical strip \(S\subset R\) whose image \(f^n(S)\)
is a horizontal strip in \(R\). The return branch \(f^n|_S\) maps one vertical side of \(S\)
into the left side of \(R\) and the other into the right side.
If two distinct source strips of order \(n\) meet, they are nested by
\cite[Lemma~11.11]{CLP}, say \(S\subsetneq S'\). Then \(f^n(V)\) meets
\(\operatorname{int}f^n(S')\) for some vertical side \(V\) of \(S\),
contradicting \(f^n(V)\subset\partial R\). Thus the strips are pairwise
disjoint. Their nonempty interiors contain distinct elements of a fixed
countable basis for \(M\), so there are at most countably many strips.

For (2), let \(\gamma\) cross \(R\) from left to right. Take the first parameter \(b\)
with \(\gamma(b)\) on the right side of \(S\), then the last \(a<b\) with
\(\gamma(a)\) on its left side. Since \(S\) extends from top to bottom,
\(\gamma_0=\gamma|_{[a,b]}\) stays in \(S\), and \(f^n(\gamma_0)\) joins the vertical sides of \(R\)
(Figure~\ref{fig:return-strip}).

For (3), iterating the two-branch construction on \(S_1\cap f^{-n_1}S_2\)
gives the stated source strip \(S\).
The concatenation property gives hyperbolicity constants \(C,\kappa\) independent of \(m\)
\cite[Definitions~4.12, 4.14 and Theorem~D]{CLP}.
Admissible arcs in \(R\) are graphs of bounded slope in the chart of
\cite[(5.11) and Remark~10.2]{CLP}. Let \(L_R\) bound their lengths.
Each stable side of \(S\) maps under \(f^n\) into a stable side of \(R\),
so its image at time \(n\) is tangent to the stable cone.
The forward estimate in \cite[(4.3)]{CLP}, applied from time \(0\), therefore bounds its length
at time \(k\) by \(CL_Re^{-\kappa k}\).
Each unstable side of \(S\) lies on an unstable side of \(R\) at time \(0\),
so it is tangent to the unstable cone there. The backward estimate from
time \(n\) bounds its length at time \(k\) by \(CL_Re^{-\kappa(n-k)}\). Thus the
perimeter \(P_k\) of \(f^k(S)\) satisfies
\[
  P_k\le 2CL_R\bigl(e^{-\kappa k}+e^{-\kappa(n-k)}\bigr).
\]
To bound the distance from \(f^kx\in\operatorname{int}f^k(S)\) to the boundary,
take the unstable coordinate segment \(J_x\) from \(x\) to its first hit
\(p_x\in\partial S\).
By \cite[(4.2)--(4.3)]{CLP}, \(f^n(J_x)\) is an unstable admissible arc and
\[
  \operatorname{length}(f^kJ_x)
  \le C e^{-\kappa(n-k)}\operatorname{length}(f^nJ_x)
  \le CL_R e^{-\kappa(n-k)}.
\]
Join \(f^kx\) to \(f^ky\) along \(f^kJ_x\), the shorter boundary arc from
\(f^kp_x\) to \(f^kp_y\) (of length at most \(P_k/2\)), and \(f^kJ_y\) in reverse.
By density of the interior,
\(\diam f^k(S)\le 2CL_R e^{-\kappa(n-k)}+P_k/2\), giving \eqref{eq:thin}
with \(B=3CL_R\).
\end{proof}

\subsection*{A horseshoe in the return rectangle}
In Section~\ref{sec:passages} we use the symbolic coding of this horseshoe
to construct a quadrilateral whose sides lie in the stable and unstable
manifolds of a single periodic orbit.

\begin{figure}[htbp]
\centering
\begin{tikzpicture}[x=.88cm,y=.88cm,font=\footnotesize]
  \coordinate (saddle) at (2.45,1.12);
  \fill[black!9] (1.70,0) rectangle (3.22,2.30);
  \fill[black!9] (0,.70) rectangle (5.65,1.49);
  \draw[black!50] (1.70,0) -- (1.70,2.30) (3.22,0) -- (3.22,2.30)
                  (0,.70) -- (5.65,.70) (0,1.49) -- (5.65,1.49);
  \draw (0,0) rectangle (5.65,2.30);
  \draw[densely dashed,line width=.85pt] (2.32,0)
    .. controls (2.34,.40) and (2.48,.80) .. (saddle)
    .. controls (2.42,1.48) and (2.24,1.91) .. (2.29,2.30);
  \draw[line width=.85pt] (0,1.03)
    .. controls (.85,.92) and (1.76,1.07) .. (saddle)
    .. controls (3.32,1.18) and (4.39,1.29) .. (5.65,1.17);
  \foreach \pt in {(2.32,0),(2.29,2.30),(0,1.03),(5.65,1.17)}
    \fill \pt circle (1pt);
  \fill (saddle) circle (1.6pt) node[below right,inner sep=2pt] {\(z\)};
  \node at (2.79,1.96) {\(S_0\)};
  \node at (4.66,.93) {\(F(S_0)\)};
  \node at (5.30,1.93) {\(R\)};
\end{tikzpicture}
\caption{The stable arc in \(S_0\) (dashed) and the unstable arc in
\(F(S_0)\) (solid) meet at the \(F\)-fixed saddle \(z\).}
\label{fig:saddle-arcs}
\end{figure}
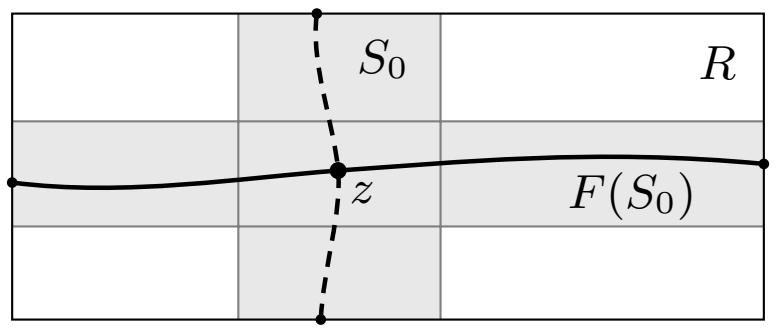

\begin{lemma}\label{lem:anchor}
Let \(R\) be a return rectangle supplied by Lemma~\ref{lem:branches}.
For some \(n\ge1\), two
return branches of order \(n\) in \(R\) generate a horseshoe for \(f^n\)
containing a fixed saddle \(z\in\operatorname{int}R\).
The saddle \(z\) has full-length stable and unstable arcs in \(R\).
These are compact arcs through \(z\) in \(W^s(z)\) and \(W^u(z)\)
joining the two horizontal and the two vertical sides of \(R\), respectively.
At each endpoint the arc meets the relative interior of
the corresponding side transversely.
\end{lemma}
\begin{proof}
Let \(\mu,A,T,B\) be as in Lemma~\ref{lem:branches} for this rectangle.
Since \(\mu\) is nonatomic, \(\mu^{\otimes6}\)-almost every tuple in
\(A^6\) has distinct coordinates. By recurrence for \((f^T)^{\times6}\),
choose such a tuple \((x_1,\ldots,x_6)\) returning to \(A^6\) in both
directions. Let \(C,\kappa\) be the uniform branch constants from the proof of
Lemma~\ref{lem:branches}. Choose \(a,b\in T\mathbb N\) with
\(f^{-a}x_i,f^bx_i\in A\) for all \(i\), taking both return times large enough that
\[
 B(e^{-\kappa a}+e^{-\kappa b})<\min_{i<j}d(x_i,x_j),
 \qquad C^{-1}e^{\kappa(a+b)}>1.
\]
Set \(n=a+b\). Since \(f^{-a}x_i\in A\cap f^{-n}A\),
Lemma~\ref{lem:branches}(1) gives a source strip of order \(n\) through each
of these six points. By \eqref{eq:thin} at \(k=a\) and the choice of \(a,b\),
no such strip contains two of them.
Only the leftmost and rightmost can meet a vertical side of \(R\).
Likewise, at most two of their images meet a horizontal side.
We may therefore choose two remaining strips \(S_0,S_1\) and put \(F=f^n\).

Put \(U=S_0\cup S_1\) and \(\Lambda=\bigcap_{j\in\mathbb Z}F^{-j}U\).
For a word \(\omega=(\omega_{-m},\ldots,\omega_{m-1})\), the nonempty compact source strip
at time \(-m\) is \(S(\omega)=\bigcap_{r=0}^{2m-1}F^{-r}S_{\omega_{-m+r}}\).
The \(2m\) branches span \(2mn\) iterates of \(f\). From time \(-m\) to time \(0\)
we apply \(F^m=f^{mn}\), so \eqref{eq:thin} with \(k=mn\) gives
\(\diam F^m(S(\omega))\le 2Be^{-\kappa mn}\).
For each bi-infinite word the
nested sets therefore determine a unique point of \(\Lambda\),
continuously in the word. Disjointness of \(S_0,S_1\) gives injectivity,
so the coding conjugates the full two-shift to \(F|_\Lambda\).
By \cite[(4.2)--(4.3)]{CLP}, \(DF\) expands the unstable cone and
\(DF^{-1}\) expands the stable cone, each by a factor of at least \(C^{-1}e^{\kappa n}>1\),
giving uniform hyperbolicity. Moreover, \(\Lambda\subset U\cap F(U)\subset
\operatorname{int}R\). Each side of \(S_0,S_1\) either lies on
\(\partial R\) or maps there, so \(\operatorname{int}U\) isolates \(\Lambda\).

The all-zero itinerary gives an \(F\)-fixed saddle \(z\).
The full-length stable arc \(\bigcap_{k\ge0}F^{-k}(S_0)\) in \(S_0\)
and the full-length unstable arc \(\bigcap_{k\ge1}F^k(S_0)\) in \(F(S_0)\)
from \cite[Proposition~4.10]{CLP} both contain \(z\).
As \(k\to\infty\), \(F^kx\to z\) on the first arc and \(F^{-k}x\to z\) on
the second. Thus they lie in \(W^s(z)\) and \(W^u(z)\), respectively
(Figure~\ref{fig:saddle-arcs}). The endpoints avoid the corners by construction.
At each endpoint \(x\) of the stable arc, a one-sided difference quotient
in \cite[(1.4)]{CLP} gives \(\|D_xFv\|\le Ce^{-\kappa n}\|v\|\) on its tangent line.
The horizontal side through \(x\) is tangent to the unstable cone, so
\(\|D_xFv\|\ge C^{-1}e^{\kappa n}\|v\|\) on its tangent line by \cite[(4.3)]{CLP}.
By our choice of \(n=a+b\), we have \(Ce^{-\kappa n}<1\), so these estimates
for a single return \(F=f^n\) rule out tangency.
At the endpoints of the unstable arc in \(F(S_0)\), the same argument for
\(F^{-1}\) gives transversality to the vertical sides.
\end{proof}

\FloatBarrier
\section{Passages of fixed duration between rectangles}\label{sec:passages}

Let \(R,R'\) be closed rectangles with designated vertical sides. For a
positive integer \(t\), we say there is a \emph{passage} of duration \(t\)
from \(R\) to \(R'\), written \(R\leadsto_t R'\), if every horizontal
crossing of \(R\) has a subpath \(\eta\) such that \(f^t\circ\eta\) is a
horizontal crossing of \(R'\). Choosing subpaths successively gives
\[
  R\leadsto_s R' \quad\text{and}\quad R'\leadsto_t R''
  \quad\Longrightarrow\quad R\leadsto_{s+t}R''.
\]

\begin{proposition}\label{prop:passage}
Let \(R_1,R_2\) be return rectangles. Suppose that in each rectangle
two return branches generate a horseshoe containing a periodic saddle
as in Lemma~\ref{lem:anchor}.
If \(f\) is topologically transitive, then there exist positive integers
\(u,v\) such that
\[
  R_1\leadsto_u R_2 \qquad\text{and}\qquad R_2\leadsto_v R_1.
\]
\end{proposition}

\begin{proof}
In each \(R_i\), the two return branches generating the horseshoe are restrictions
of the same iterate \(f^{n_i}\). Let \(O_i\) be the \(f\)-orbit of the all-zero
saddle. Fix a horseshoe point \(x_i\notin O_1\cup O_2\).
For \(\ell\in\mathbb N\), the past block of its \(f^{n_i}\)-itinerary on \([-\ell,-1]\)
has source strip \(P_\ell\) at time \(-\ell\) and horizontal image strip
\(H_\ell^-=f^{\ell n_i}(P_\ell)\) at time \(0\). The future block \(w\) on
\([0,\ell-1]\) has vertical source strip \(V_\ell^+\). The concatenated
branch therefore has source strip \(S_\ell\) given by
\[
 \begin{aligned}
 S_\ell &= P_\ell\cap f^{-\ell n_i}(V_\ell^+),\\
 f^{\ell n_i}(S_\ell) &= f^{\ell n_i}(P_\ell)\cap V_\ell^+
                        = H_\ell^-\cap V_\ell^+.
 \end{aligned}
\]
Applying \eqref{eq:thin} at the midpoint of this branch of order \(2\ell n_i\)
gives \(\diam(H_\ell^-\cap V_\ell^+)\le 2Be^{-\kappa\ell n_i}\).
Since \(x_i\in H_\ell^-\cap V_\ell^+\) has positive distance from
\(\partial R_i\cup O_1\cup O_2\), fix \(\ell\) large enough that
\(H_\ell^-\cap V_\ell^+\subset\operatorname{int}R_i\setminus(O_1\cup O_2)\)
for both \(i=1,2\).
For each of \(w0^\infty\) and \(w10^\infty\),
the nested source strips give a full-length stable arc in \(V_\ell^+\)
by \cite[Proposition~4.10]{CLP}. After \(\ell\) iterates of \(f^{n_i}\),
their images lie in disjoint source strips. After one further iterate, both lie on the
stable arc of the all-zero saddle. Thus the arcs are disjoint and lie in
\(W^s(O_i)\). Similarly, choose two distinct extensions \(\eta_0,\eta_1\)
of the past block with symbol \(0\) at all sufficiently negative times.
They give two disjoint full-length unstable arcs in \(H_\ell^-\cap W^u(O_i)\).
Pairing \(\eta_a\) with \(w0^\infty\) and \(w10^\infty\) codes the unique
intersections \(p_{a0}\) and \(p_{a1}\), respectively. On each of the four
arcs, take the subarc between its two intersection points. These subarcs
bound an open quadrilateral \(Q_i\) with
\(\overline{Q_i}\subset H_\ell^-\cap V_\ell^+\) (Figure~\ref{fig:four-arcs}).

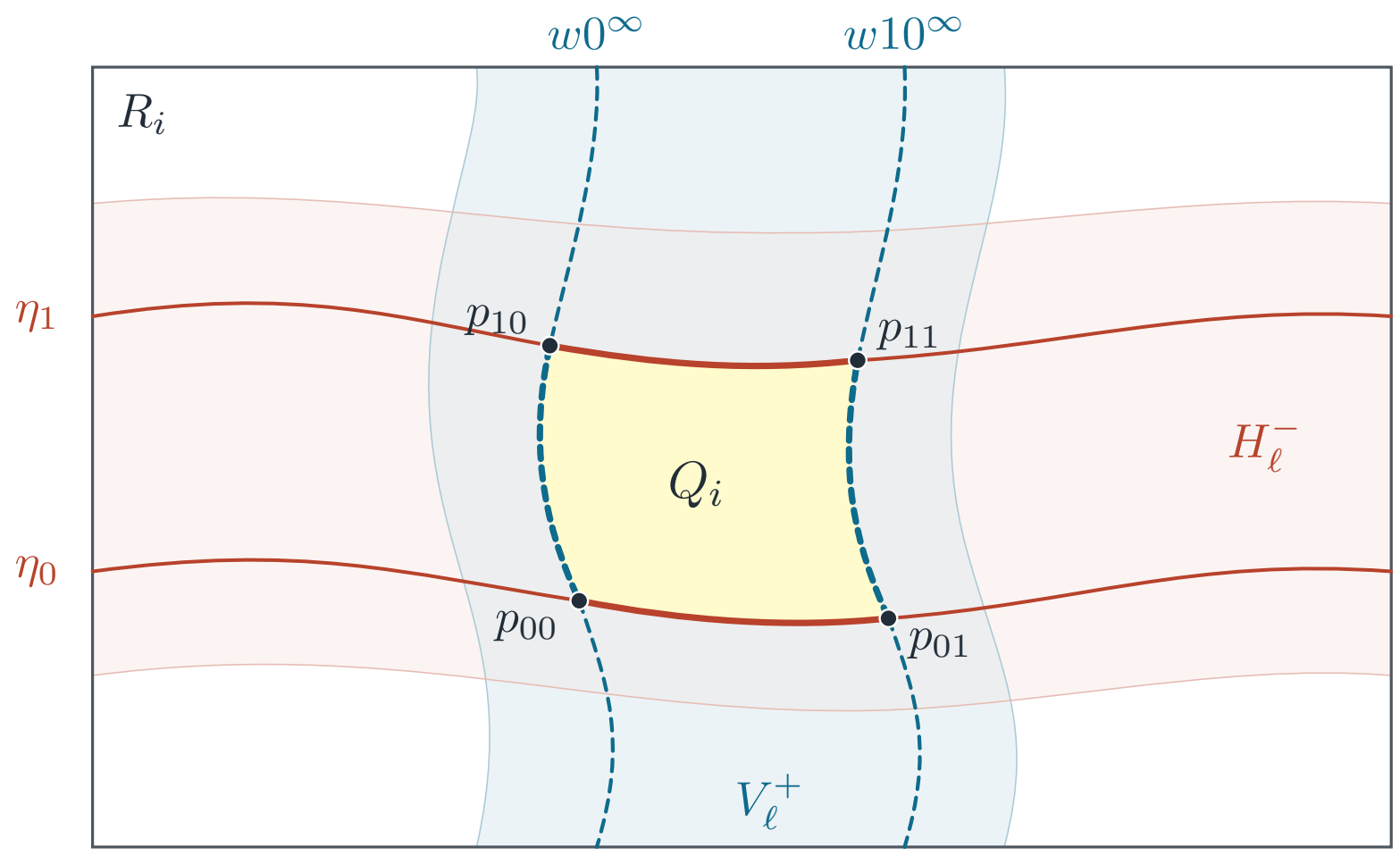
\begin{figure}[htbp]
\centering
\begin{tikzpicture}[x=.08577\linewidth,y=.08577\linewidth,font=\small,
    stable/.style={draw=stableblue,dash pattern=on 2.5pt off 1.8pt,line cap=round},
    unstable/.style={draw=unstablered,line cap=round}]
  \definecolor{stableblue}{RGB}{14,107,140}
  \definecolor{unstablered}{RGB}{184,66,43}
  \definecolor{cornerink}{RGB}{33,45,56}
  \path[use as bounding box] (-.82,-.24) rectangle (9.09,6.01);
  \def\verticalstrip{(2.62,0)
    .. controls (3.02,1.78) and (1.91,2.10) .. (2.44,4.26)
    .. controls (2.55,4.76) and (2.65,5.04) .. (2.62,5.32)
    -- (6.22,5.32)
    .. controls (6.31,4.25) and (5.60,3.38) .. (5.96,2.04)
    .. controls (6.19,1.20) and (6.44,.84) .. (6.22,0) -- cycle}
  \def\horizontalstrip{(0,1.17)
    .. controls (2.00,1.43) and (3.33,.94) .. (5.07,.93)
    .. controls (6.28,.91) and (7.65,1.27) .. (8.86,1.17)
    -- (8.86,4.39)
    .. controls (7.47,4.47) and (6.11,4.18) .. (4.78,4.19)
    .. controls (2.88,4.17) and (1.80,4.55) .. (0,4.39) -- cycle}
  \fill[stableblue!8] \verticalstrip;
  \fill[unstablered!6] \horizontalstrip;
  \begin{scope}
    \clip \verticalstrip;
    \fill[cornerink!8] \horizontalstrip;
  \end{scope}
  \draw[stableblue!35,line width=.3pt] \verticalstrip;
  \draw[unstablered!35,line width=.3pt] \horizontalstrip;
  \draw[cornerink!80,line width=.6pt] (0,0) rectangle (8.86,5.32);

  \coordinate (p00) at (3.32,1.68);
  \coordinate (p01) at (5.43,1.56);
  \coordinate (p10) at (3.12,3.42);
  \coordinate (p11) at (5.22,3.32);
  \def\bottomside{.. controls (4.04,1.54) and (4.77,1.49) .. (p01)}
  \def\rightside{.. controls (5.20,2.06) and (5.08,2.57) .. (p11)}
  \def\topside{.. controls (4.49,3.24) and (3.87,3.29) .. (p10)}
  \def\leftside{.. controls (2.98,2.78) and (3.06,2.24) .. (p00)}
  \fill[yellow!25] (p00) \bottomside \rightside \topside \leftside -- cycle;

  \draw[stable,line width=.8pt] (3.44,5.32)
    .. controls (3.47,4.71) and (3.27,4.07) .. (p10);
  \draw[stable,line width=.8pt] (5.54,5.32)
    .. controls (5.57,4.62) and (5.36,3.99) .. (p11);
  \draw[stable,line width=.8pt] (p00)
    .. controls (3.51,1.12) and (3.66,.73) .. (3.44,0);
  \draw[stable,line width=.8pt] (p01)
    .. controls (5.64,.97) and (5.73,.65) .. (5.54,0);
  \draw[unstable,line width=.75pt] (0,1.88)
    .. controls (1.44,2.07) and (2.06,1.88) .. (p00)
    (p01) .. controls (6.73,1.66) and (7.49,1.98) .. (8.86,1.88);
  \draw[unstable,line width=.75pt] (0,3.62)
    .. controls (1.40,3.84) and (2.09,3.61) .. (p10)
    (p11) .. controls (6.59,3.40) and (7.55,3.71) .. (8.86,3.62);
  \draw[stable,line width=1.35pt] (p10) \leftside (p01) \rightside;
  \draw[unstable,line width=1.35pt] (p00) \bottomside (p11) \topside;
  \foreach \pt/\pos/\idx in {p00/below left/00,p01/below right/01,
      p10/above left/10,p11/above right/11}
    \filldraw[fill=cornerink,draw=white,line width=.45pt] (\pt) circle (1.85pt)
      node[\pos,text=cornerink,inner xsep=4pt,inner ysep=2pt] {\(p_{\idx}\)};
  \node[above,text=stableblue] at (3.44,5.32) {\(w0^\infty\)};
  \node[above,text=stableblue] at (5.54,5.32) {\(w10^\infty\)};
  \node[left,text=unstablered] at (-.10,1.88) {\(\eta_0\)};
  \node[left,text=unstablered] at (-.10,3.62) {\(\eta_1\)};
  \node[text=cornerink] at (.34,5.00) {\(R_i\)};
  \node[text=unstablered] at (8.00,2.75) {\(H_\ell^-\)};
  \node[text=stableblue] at (4.62,.31) {\(V_\ell^+\)};
  \node[text=cornerink,font=\normalsize] at (4.12,2.47) {\(Q_i\)};
\end{tikzpicture}
\caption{The dashed stable arcs correspond to the futures \(w0^\infty\) and
\(w10^\infty\) from left to right. The solid unstable arcs correspond to the
pasts \(\eta_0\) and \(\eta_1\) from bottom to top. The marked points
\(p_{ab}\), indexed by past \(a\) and future \(b\), are the corners of \(Q_i\).}
\label{fig:four-arcs}
\end{figure}

Write \(\partial^s Q_i,\partial^u Q_i\) for the unions of the closed stable
and unstable sides. By construction, \(\overline{Q_2}\) is at positive distance
from \(O_1\), so the forward iterates of \(\partial^s Q_1\), which approach
\(O_1\), eventually miss \(\overline{Q_2}\). Likewise, \(\overline{Q_1}\) is
at positive distance from \(O_2\), so the backward iterates of \(\partial^u Q_2\),
which approach \(O_2\), eventually miss \(\overline{Q_1}\). Thus there is
\(N\ge 1\) such that for every \(n\ge N\),
\[
 f^n(\partial^s Q_1)\cap\overline{Q_2}=\varnothing,\qquad
 f^n(\overline{Q_1})\cap\partial^u Q_2=\varnothing .
\]
Applying transitivity to the nonempty open sets \(f^N(Q_1)\) and \(Q_2\)
gives \(j\ge0\) with \(f^j(f^N(Q_1))\cap Q_2\ne\varnothing\).
Set \(n=N+j\).
As in \cite[proof of Theorem~6.4, Step~2]{BCS}, if \(f^n(\partial^u Q_1)\)
missed \(Q_2\), connectedness would give \(Q_2\subset f^n(Q_1)\), hence
\(\partial^u Q_2\subset f^n(\overline{Q_1})\), a contradiction.
Thus the image \(\alpha\) of an unstable side of \(Q_1\) meets \(Q_2\).
Its endpoints lie outside \(\overline{Q_2}\), and it misses the closed
unstable sides of \(Q_2\) (Figure~\ref{fig:quadrilateral-entry}).

\begin{figure}[htbp]
\centering
\begin{tikzpicture}[x=1cm,y=1cm,font=\footnotesize,
    stable/.style={densely dashed},unstable/.style={line width=.55pt},
    widened/.style={draw=black!65,densely dotted,line width=.7pt}]
  \definecolor{entryblue}{RGB}{24,99,133}
  \def\outerarc{(.1,2.05) .. controls (5.7,2.05) and (5.7,.25) .. (.1,.25)}
  \def\innerarc{(.1,.78) .. controls (3.85,.78) and (3.85,1.52) .. (.1,1.52)}
  \path \outerarc coordinate[pos=.20] (entry) coordinate[pos=.40] (target)
                  coordinate[pos=.67] (lowerarc);
  \coordinate (base) at (0,0);
  \coordinate (top) at (0,2.30);
  \fill[black!12] \outerarc -- \innerarc -- cycle;
  \draw[black!60,line width=.55pt] (2.05,-.58) rectangle (8.10,3.00);
  \node at (7.75,2.72) {$R_2$};
  \draw[unstable] \outerarc;
  \draw[unstable] \innerarc;
  \draw[stable] (.1,.25) -- (.1,.78) (.1,1.52) -- (.1,2.05);
  \coordinate (wideleft) at (entry |- base);
  \path (wideleft) -- (4.975,0) coordinate[pos=2] (wideright);
  \draw[widened] (wideleft) rectangle (wideright |- top);
  \draw[unstable] (3.25,0) -- (6.70,0) (3.25,2.30) -- (6.70,2.30);
  \draw[stable] (6.70,0) -- (6.70,2.30);
  \draw[stable,line width=.9pt] (3.25,-.58) -- (3.25,3.00)
    node[pos=.91,right,inner sep=2pt] {$\ell^s$};
  \begin{scope}
    \clip (entry |- target) rectangle (6.70,2.30);
    \draw[entryblue,line width=1.6pt] \outerarc;
  \end{scope}
  \fill (entry) circle (1.35pt) node[above left,inner sep=1.5pt] {$a$};
  \fill (target) circle (1.35pt) node[right,inner sep=2.5pt] {$y$};
  \node[text=entryblue] at (3.85,1.98) {$f^n(I)$};
  \node at (1.25,1.76) {$f^n(Q_1)$};
  \node[below right,inner sep=1.5pt] at (lowerarc) {$\alpha$};
  \node at (6.24,1.94) {$Q_2$};
  \path (wideright) -- (wideright |- top)
    node[midway,right,inner sep=2pt] {$\widetilde Q_2$};
\end{tikzpicture}
\caption{The image \(f^n(Q_1)\) is shaded, and the dotted lines show the
enlargement \(\widetilde Q_2\). Here \(\alpha\) enters and leaves \(Q_2\)
through its left side. The dashed stable arc \(\ell^s\) extends this side
to the top and bottom of \(R_2\). The bold subarc \(f^n(I)\) joins
\(a\in\partial\widetilde Q_2\) and \(y\in Q_2\) in opposite components of
\(R_2\setminus\ell^s\).}
\label{fig:quadrilateral-entry}
\end{figure}
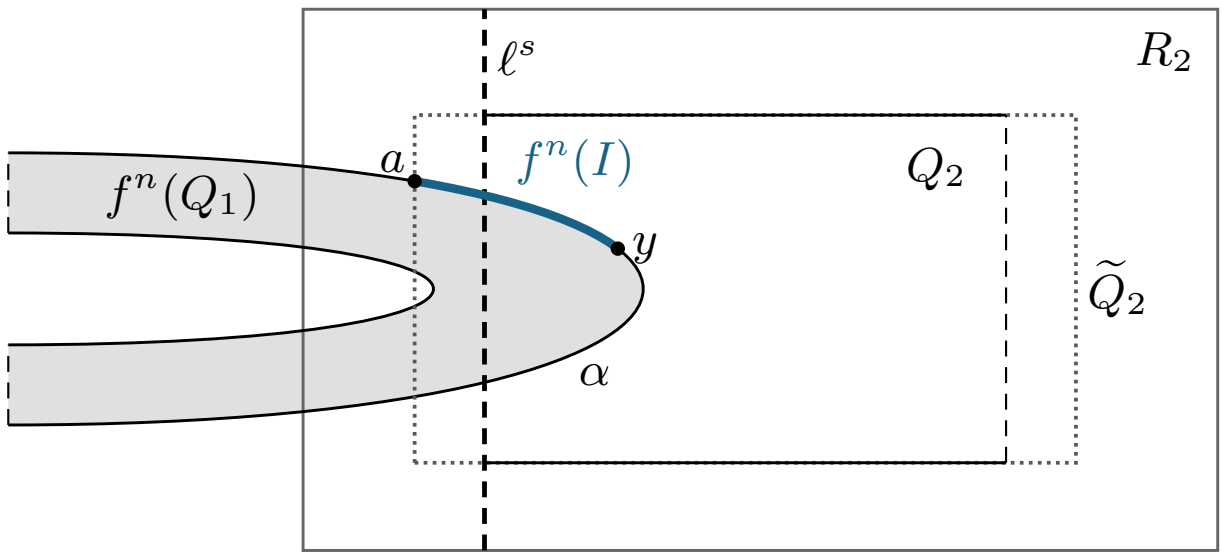

Use rectangle coordinates near \(\overline{Q_2}\) in which
\(Q_2=(0,1)^2\) and its two stable sides extend along \(x=0,1\).
Set \(\widetilde Q_2=(-\varepsilon,1+\varepsilon)\times(0,1)\).
By compactness, for sufficiently small \(\varepsilon>0\),
\(\overline{\widetilde Q_2}\subset\operatorname{int}R_2\), the endpoints of
\(\alpha\) lie outside \(\overline{\widetilde Q_2}\), and \(\alpha\) misses
both horizontal edges of \(\widetilde Q_2\). Follow \(\alpha\) from an endpoint
to \(y\in\alpha\cap Q_2\).
Let \(a\) be the last point where \(\alpha\) meets \(\partial\widetilde Q_2\)
before reaching \(y\). The subarc from \(a\) to \(y\) lies in
\(\overline{\widetilde Q_2}\) and has the form \(f^n(I)\), with \(I\) on the
unstable side of \(Q_1\). Since \(a\) lies on \(x=-\varepsilon\) or
\(x=1+\varepsilon\), take \(\ell^s\) to be the full-length
stable arc in \(R_2\) containing the corresponding side of \(Q_2\).
This arc joins the two horizontal sides of \(R_2\) and separates \(R_2\)
into two components. The compact arc \(f^n(I)\subset\operatorname{int}R_2\)
has endpoints \(a\) and \(y\) in opposite components of \(R_2\setminus\ell^s\)
(Figure~\ref{fig:quadrilateral-entry}).

We now choose finite return strips so that every path in a thin strip
along \(I\) joining its two end sections has an \(f^n\)-image crossing a
fixed vertical strip in \(R_2\). Let \(S_0,S_1\) be the source strips of the two horseshoe branches of order
\(n_2\) in \(R_2\), and let \((\omega_j)_{j\ge0}\) be the future itinerary
defining \(\ell^s\). The source strip for its first \(k\) symbols is
\[
  V_k=\bigcap_{j=0}^{k-1}f^{-j n_2}(S_{\omega_j}),\qquad k\ge1.
\]
These strips decrease to \(\ell^s\). Since \(a,y\notin\ell^s\), fix \(k\)
large enough that both lie outside \(V_k\).
They remain in different components of \(R_2\setminus V_k\), since
\(R_2\setminus V_k\subset R_2\setminus\ell^s\).

Let \(\ell^u\) be the full-length unstable arc in \(R_1\) containing the
side of \(Q_1\) on which \(I\) lies. The last \(m\) symbols of its past
itinerary give image strips \(H_m\) decreasing to \(\ell^u\).
Choose coordinates on \(R_1\) with vertical sides \(x=0,1\) in which
\(\ell^u=[0,1]\times\{1/2\}\) and \(I=[s,t]\times\{1/2\}\).
The images of the two endpoints are \(a\) and \(y\) in some order.
With \(k\) fixed, choose open disks \(U_s,U_t\) about \(f^n(s,1/2)\) and
\(f^n(t,1/2)\), respectively, in their distinct components of \(R_2\setminus V_k\).
Put \(K_m=H_m\cap\{s\le x\le t\}\). The compact sets \(f^n(K_m)\)
decrease to \(f^n(I)\subset\operatorname{int}R_2\), and for \(r=s,t\)
the sets \(f^n(K_m\cap\{x=r\})\) decrease to \(\{f^n(r,1/2)\}\).
Thus we can fix \(m\) large enough that \(f^n(K_m)\subset\operatorname{int}R_2\)
and \(f^n(K_m\cap\{x=r\})\subset U_r\) for \(r=s,t\)
(Figure~\ref{fig:stable-side}).

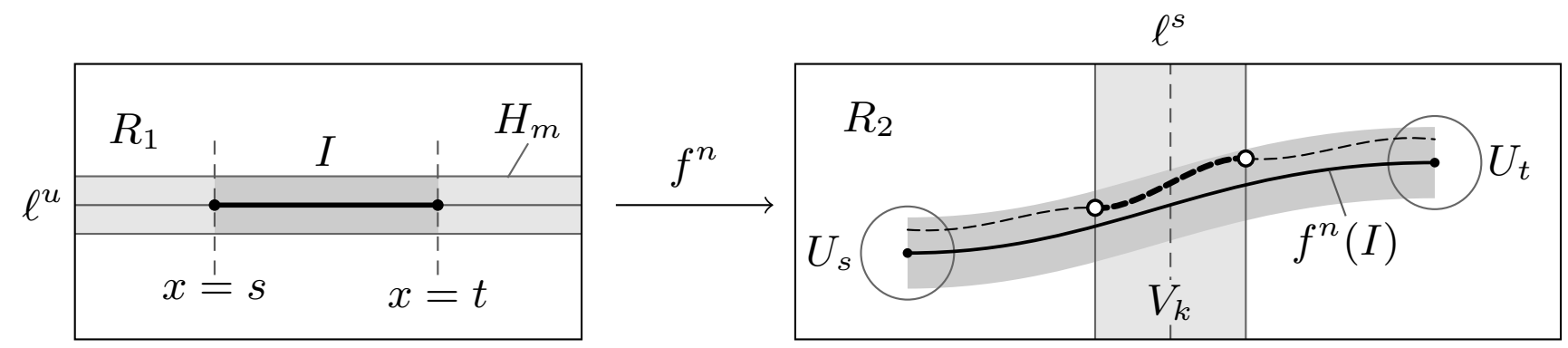
\begin{figure}[htbp]
\centering
\begin{tikzpicture}[x=1cm,y=1cm,font=\footnotesize,
    arc part/.style={line width=1.05pt,line cap=round},
    fiber/.style={densely dashed,black!65}]
  \fill[black!10] (0,.92) rectangle (3.70,1.34);
  \fill[black!20] (1.02,.92) rectangle (2.65,1.34);
  \draw[black!60] (0,.92) -- (3.70,.92) (0,1.34) -- (3.70,1.34);
  \draw[black!65] (0,1.13) -- (3.70,1.13);
  \node[left,inner sep=2pt] at (0,1.13) {$\ell^u$};
  \draw[fiber] (1.02,.62) -- (1.02,1.60) (2.65,.62) -- (2.65,1.60);
  \node[below,inner sep=1pt] at (1.02,.62) {$x=s$};
  \node[below,inner sep=1pt] at (2.65,.62) {$x=t$};
  \draw[arc part] (1.02,1.13) -- (2.65,1.13);
  \fill (1.02,1.13) circle (1.2pt) (2.65,1.13) circle (1.2pt);
  \node[above,inner sep=2pt] at (1.83,1.34) {$I$};
  \node at (.44,1.67) {$R_1$};
  \draw[black!60] (3.16,1.33) -- (3.31,1.58);
  \node[above,inner sep=1pt] at (3.31,1.58) {$H_m$};
  \draw[->] (3.95,1.13) -- node[above] {$f^n$} (5.10,1.13);
  \fill[black!10] (7.45,.15) rectangle (8.55,2.16);
  \def\imagearc{(6.08,.78) .. controls (7.55,.78) and (8.25,1.44) .. (9.93,1.44)}
  \path[fill=black!20] {[yshift=-.26cm] \imagearc}
    -- (9.93,1.70) .. controls (8.25,1.70) and (7.55,1.04) .. (6.08,1.04)
    -- cycle;
  \draw[black!60] (7.45,.15) -- (7.45,2.16) (8.55,.15) -- (8.55,2.16);
  \draw[fiber] (8,.15) -- (8,2.16);
  \node[above] at (8,2.16) {$\ell^s$};
  \node[fill=black!10,inner sep=1pt] at (8,.42) {$V_k$};
  \node at (5.80,1.76) {$R_2$};
  \draw[black!60] (6.08,.78) circle (.34) (9.93,1.44) circle (.34);
  \node[left,inner sep=1pt] at (5.74,.78) {$U_s$};
  \node[right,inner sep=1pt] at (10.27,1.44) {$U_t$};
  \draw[line width=.70pt] \imagearc;
  \fill (6.08,.78) circle (1pt) (9.93,1.44) circle (1pt);
  \draw[black!65] (9.16,1.38) -- (9.28,1.05);
  \node[below,inner sep=1pt] at (9.28,1.05) {$f^n(I)$};
  \begin{scope}[draw=black,dash pattern=on 2.5pt off 1.4pt,line cap=round]
    \draw (6.08,.95) .. controls (6.65,.91) and (7.00,1.13) .. (7.45,1.11);
    \draw[line width=1.25pt] (7.45,1.11)
      .. controls (7.78,1.08) and (8.22,1.50) .. (8.55,1.47);
    \draw (8.55,1.47) .. controls (9.04,1.43) and (9.40,1.67) .. (9.93,1.61);
  \end{scope}
  \foreach \p in {(7.45,1.11),(8.55,1.47)}
    \draw[fill=white,line width=.65pt] \p circle (1.55pt);
  \draw (0,.15) rectangle (3.70,2.16) (5.26,.15) rectangle (10.85,2.16);
\end{tikzpicture}
\caption{The bold subarc \(I\subset\ell^u\) lies in the darker strip
\(K_m\subset H_m\). On the right, the dashed path in \(f^n(K_m)\) joins \(U_s\)
to \(U_t\). The bold dashed segment between the two small open circles crosses \(V_k\)
from one vertical side to the other.}
\label{fig:stable-side}
\end{figure}

With \(m,n,k\) fixed, let \(\gamma\) be any horizontal crossing of \(R_1\).
Applying Lemma~\ref{lem:branches}(2) successively along the \(m\) prescribed
branches gives a subpath of \(\gamma\) whose \(f^{mn_1}\)-image crosses \(R_1\)
inside \(H_m\). By continuity of the \(x\)-coordinate, this crossing contains
a subpath in \(K_m\) joining \(x=s\) to \(x=t\).
By the choice of \(m\), its \(f^n\)-image lies in \(\operatorname{int}R_2\)
and joins \(U_s\) to \(U_t\). These disks lie in opposite components of
\(R_2\setminus V_k\), so the image contains a subpath in \(V_k\) joining its
two vertical sides
(Figure~\ref{fig:stable-side}).
The concatenated branch \(f^{kn_2}|_{V_k}\) of Lemma~\ref{lem:branches}(3)
maps the two vertical sides of \(V_k\) into opposite vertical sides of \(R_2\),
so it sends this subpath to a horizontal crossing of \(R_2\).
This gives \(R_1\leadsto_u R_2\)
with \(u=mn_1+n+kn_2\), independently of \(\gamma\).
Interchanging \(R_1,R_2\) gives \(v\).
\end{proof}

\section{Counting synchronized returns}\label{sec:count}

Suppose that distinct ergodic measures \(\mu_1,\mu_2\) attain
\(h=\htop(f)>0\). Each \(\mu_i\) is nonatomic since \(h>0\), and is of
saddle type by Ruelle's inequality for \(f\) and \(f^{-1}\)~\cite{Ruelle}.
Take \(R_i,A_i,T_i\) from Lemma~\ref{lem:branches}
for \(\mu_i\), using common constants \(B,\kappa\) in \eqref{eq:thin}.
By Lemma~\ref{lem:anchor} and Proposition~\ref{prop:passage},
\(R_1\leadsto_u R_2\) and \(R_2\leadsto_v R_1\) for fixed \(u,v\).
For a Lipschitz function \(\varphi\), write
\[
 b_i=\int\varphi\,d\mu_i,\qquad
 \varphi_q(z)=\frac{1}{q}\sum_{k=0}^{q-1}\varphi(f^kz).
\]
Choose \(\varphi\) so that \(b_1<b_2\), and fix \(0<\eta<(b_2-b_1)/4\).

\begin{samepage}
Summing \eqref{eq:thin} and using the Lipschitz bound gives a uniform
constant \(C_\varphi\) such that each return branch \(f^n|_S\)
in either rectangle satisfies
\[
 \sum_{k=0}^{n-1}|\varphi(f^kz)-\varphi(f^kx)|\le C_\varphi
 \qquad (x,z\in S).
\]
\end{samepage}
\begin{samepage}
By Birkhoff's ergodic theorem applied to \(f^{-1}\) and \(f\), respectively,
we can choose \(K>C_\varphi/\eta\) and sets \(D_i\subset A_i\) with
\(\mu_i(D_i)>0\) such that for all \(q\ge K\),
\[
 |\varphi_q(f^{-q}x)-b_1|<\eta \quad(x\in D_1),\qquad
 |\varphi_q(y)-b_2|<\eta \quad(y\in D_2).
\]
\end{samepage}

\subsection{A lower bound for the entropy of return labels}

Fix a common multiple \(Q\) of \(T_1,T_2\) with \(Q-u\ge K\).
Set \(D=D_1\times D_2\) and \(\rho=\mu_1\otimes\mu_2\).
For each multiple \(L\ge 2Q\) of \(Q\), consider all ways to split \(L\)
into an ordered pair \((j,L-j)\) of positive multiples of \(Q\).
The possible values of \(j\) form \(J_L=\{Q,2Q,\ldots,L-Q\}\), with
\(\lvert J_L\rvert=L/Q-1\). For each \(j\in J_L\), the first return time
to \(D\) under \(G_j=f^j\times f^{L-j}\) is given almost everywhere by
\[
 \tau_j(x,y)=\min\{m\ge1:f^{mj}x\in D_1,\ f^{m(L-j)}y\in D_2\}.
\]
We use the following lemma to compute the entropy of the induced map.
\begin{lemma}\label{lem:induce}
Let \(G\) be an invertible measure-preserving transformation of a standard
probability space \((X,\mathcal B,\rho)\). For a measurable set \(D\) with
\(\rho(D)>0\), let \(\tau\) be the first positive return time to \(D\),
write \(G_D(x)=G^{\tau(x)}x\), and put \(\rho_D=\rho|_D/\rho(D)\).
Then \(\int_D \tau\,d\rho_D\le 1/\rho(D)\).
If \(h_\nu(G)=b<\infty\) for almost every ergodic component \(\nu\) of
\(\rho\), then
\[
  h_{\rho_D}(G_D)=b\int_D \tau\,d\rho_D.
\]
\end{lemma}

\begin{proof}
Let \(Y=\bigcup_{k\in\mathbb Z}G^kD\) and \(\rho_Y=\rho|_Y/\rho(Y)\).
By recurrence, the forward iterates of \(D\) cover \(Y\) modulo \(\rho_Y\).
Applying Kac's formula~\cite[Theorem~1.2.2]{VianaOliveira} on \(Y\) gives
\[
  \int_D \tau\,d\rho_D=\frac{1}{\rho_Y(D)}=\frac{\rho(Y)}{\rho(D)}.
\]
Under the entropy hypothesis, \(h_{\rho_Y}(G|_Y)=b\) by invariance of
\(Y\) and ergodic decomposition. Abramov's formula on \(Y\)
\cite[Theorem~2]{Abramov} now gives
\(\rho_Y(D)\,h_{\rho_D}(G_D)=b\), proving the second assertion.
\end{proof}

By the power and product formulas, \(h_\rho(G_j)=\htop(G_j)=Lh\).
Hence \(h_\nu(G_j)=Lh\) for almost every \(\nu\) in the ergodic
decomposition of \(\rho\). Lemma~\ref{lem:induce} therefore gives
\begin{equation}\label{eq:synchronous-entropy}
 h_{\rho_D}((G_j)_D)=Lh\int_D \tau_j\,d\rho_D,\qquad
 \int_D \tau_j\,d\rho_D\le 1/\rho(D).
\end{equation}
For each \(m\ge1\), let \(S\) and \(T\) range over the source strips of
orders \(mj\) in \(R_1\) and \(m(L-j)\) in \(R_2\), respectively. Define
\[
  \mathcal P_j
  =\Bigl\{\,\{(x,y)\in D:\tau_j(x,y)=m,\ x\in S,\ y\in T\}\,\Bigr\}_{m,S,T},
\]
omitting sets of \(\rho_D\)-measure zero. Since \(T_1\mid mj\) and
\(T_2\mid m(L-j)\), Lemma~\ref{lem:branches}(1) gives each returning pair
a unique label \((m,S,T)\). Thus \(\mathcal P_j\) is a countable partition
of \(D\) modulo \(\rho_D\).

By Poincar\'e recurrence for \(G_j\) and \(G_j^{-1}\), the two-sided
\(\mathcal P_j\)-itinerary is defined on a set of full \(\rho_D\)-measure.
Suppose two points in this set have the same itinerary.
In either coordinate, let \(V\) be the source strip for their common block
from time \(-a\) to \(b\) under \(f\). At time zero both orbits in that
coordinate lie in \(f^a(V)\), and \eqref{eq:thin} gives
\[
 \diam f^a(V)\le B\bigl(e^{-\kappa a}+e^{-\kappa b}\bigr).
\]
Letting \(a,b\to\infty\) identifies both coordinates. Thus the two-sided
itineraries separate points on this full-measure set, so \(\mathcal P_j\)
generates \((G_j)_D\) modulo \(\rho_D\) by \cite[Proposition~9.6]{DGS}.
If the Shannon entropy \(H_{\rho_D}(\mathcal P_j)\) is finite,
subadditivity, the Kolmogorov--Sinai
theorem~\cite[Proposition~10.9(a) and Theorem~11.12]{DGS}, and
\eqref{eq:synchronous-entropy} give, in turn,
\begin{equation}\label{eq:label-lower}
 H_{\rho_D}(\mathcal P_j)
 \ge h_{\rho_D}((G_j)_D,\mathcal P_j)
 = h_{\rho_D}((G_j)_D)
 =Lh\int_D \tau_j\,d\rho_D.
\end{equation}
The resulting lower bound on \(H_{\rho_D}(\mathcal P_j)\) also holds
when this entropy is infinite.

\subsection{Distinguishing switching times and counting labels}

For each fixed split \(j\), \eqref{eq:label-lower} gives a lower bound for
the entropy of the return labels of \(G_j\). We now fix the return time
\(m\ge1\) and count labels across all splits \(j\in J_L\) by realizing them
as \(f\)-orbit segments of the common length \(N=mL+u+v\).
Let \(N_m\) count the labels \((j,m,S,T)\) with \(j\in J_L\) and
\(\rho_D(\{\tau_j=m\}\cap(S\times T))>0\). For each such label, put
\[
 F_{j,S,T}=
 S\cap f^{-(mj+u)}T\cap f^{-N}R_1 .
\]
For each \(z\in F_{j,S,T}\), the successive iterates shown below take \(z\)
through the indicated rectangles.
\[
 R_1 \xrightarrow{\ f^{mj}|_S\ } R_1 \xrightarrow{\ u\ } R_2
     \xrightarrow{\ f^{m(L-j)}|_T\ } R_2 \xrightarrow{\ v\ } R_1 .
\]
The passage from \(R_1\) to \(R_2\) begins at time \(mj\).
The second branch starts at \(mj+u\) and ends at \(mL+u\) for every \(j\).
Successively restricting a horizontal crossing of \(R_1\) along these
four steps gives a subpath in \(F_{j,S,T}\) whose \(f^N\)-image crosses
\(R_1\) horizontally.

For fixed \(j\), distinct labels have disjoint strips \(S,S'\) or disjoint
strips \(T,T'\), so the corresponding compact sets \(F_{j,S,T}\) are disjoint.

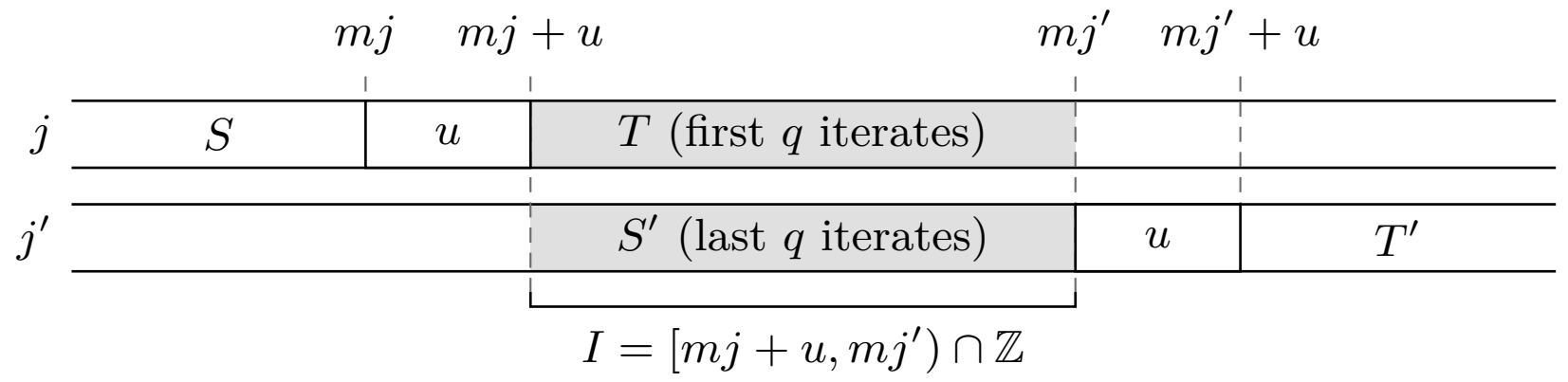
\begin{figure}[htbp]
\centering
\begin{tikzpicture}[x=1cm,y=.85cm,font=\footnotesize,
    time guide/.style={densely dashed,black!58,line width=.4pt}]
  \foreach \x/\tick in {1.65/{mj},2.8/{mj+u},6.6/{mj'},7.75/{mj'+u}} {
    \node[above] at (\x,1.62) {$\tick$};
  }
  \foreach \y in {0,.85} {
    \fill[black!12] (2.8,\y) rectangle (6.6,\y+.55);
    \draw[line width=.5pt] (-.40,\y) -- (9.95,\y)
                           (-.40,\y+.55) -- (9.95,\y+.55);
  }
  \draw[time guide] (1.65,1.60) -- (1.65,1.40);
  \draw[time guide] (7.75,1.60) -- (7.75,.55);
  \foreach \x in {2.8,6.6}
    \draw[time guide] (\x,1.60) -- (\x,-.17);
  \draw[fill=white,line width=.5pt] (1.65,.85) rectangle (2.8,1.40);
  \draw[fill=white,line width=.5pt] (6.6,0) rectangle (7.75,.55);
  \node[left] at (-.40,1.125) {$j$};
  \node[left] at (-.40,.275) {$j'$};
  \node at (.625,1.125) {\(S\)};
  \node at (2.225,1.125) {$u$};
  \node at (4.7,1.125) {\(T\) (first \(q\) iterates)};
  \node at (4.7,.275) {\(S'\) (last \(q\) iterates)};
  \node at (7.175,.275) {$u$};
  \node at (8.85,.275) {\(T'\)};
  \draw[line width=.5pt] (2.8,-.17) -- (2.8,-.29) -- (6.6,-.29) -- (6.6,-.17);
  \node[below] at (4.7,-.29) {$I=[mj+u,mj')\cap\mathbb Z$};
\end{tikzpicture}
\caption{Each branch is labeled by its source strip.
The boxes marked \(u\) show the passages from \(R_1\) to \(R_2\).
If one orbit realized both itineraries, its average over the
shaded interval \(I\) would be within \(2\eta\) of \(b_2\) in the upper row
and of \(b_1\) in the lower row.}
\label{fig:switching-times}
\end{figure}

\begin{samepage}
Suppose \(j<j'\) and \(z\in F_{j,S,T}\cap F_{j',S',T'}\).
After the passage of duration \(u\), the spacing \(j'-j\ge Q\) leaves
an overlap of at least \(K\) iterates. Its interval \(I\) and length \(q\) satisfy
\[
 \begin{gathered}
 I=[mj+u,mj')\cap\mathbb Z,\qquad q=|I|=m(j'-j)-u,\\
 K\le Q-u\le q\le\min\{m(L-j),mj'\}.
 \end{gathered}
\]
\end{samepage}
Along the orbit of \(z\), these times index the first \(q\) iterates of
the branch with source \(T\) and the last \(q\) iterates of the branch with source \(S'\)
(Figure~\ref{fig:switching-times}). Take \(y\in T\cap D_2\) as the second
coordinate of a point in the \(\mathcal P_j\)-atom labeled \((m,S,T)\),
and \(x\in S'\) as the first coordinate of a point in the
\(\mathcal P_{j'}\)-atom labeled \((m,S',T')\). The latter point returns
to \(D_1\times D_2\) after \(m\) iterates of \(G_{j'}\), so its first
coordinate at return, \(f^{mj'}x\), lies in \(D_1\).
Let \(a=\varphi_q(f^{mj+u}z)\) be the orbit average over \(I\).
Since \(mj+u=mj'-q\), applying the branch estimate on \(T\) and \(S'\),
respectively, gives
\[
 |a-\varphi_q(y)|\le C_\varphi/q,\qquad
 |a-\varphi_q(f^{mj'-q}x)|\le C_\varphi/q.
\]
As \(q\ge K\), the defining bounds for \(y\in D_2\) and \(f^{mj'}x\in D_1\)
yield \(|a-b_2|<2\eta\) and \(|a-b_1|<2\eta\).
This gives \(b_2-b_1<4\eta\), contrary to the choice of \(\eta\).
Pavlov likewise distinguishes gluing times by comparing initial and terminal
orbit segments~\cite[proofs of Theorems~1.2 and~1.3]{Pavlov}.

At this return time, the compact sets \(F_{j,S,T}\) are therefore pairwise
disjoint across all splits. Choose any \(d\ge2\) of them, denoted by
\(F_1,\ldots,F_d\), and let \(\delta>0\) be their minimum pairwise distance.
Fix a horizontal crossing \(\gamma:[0,1]\to R_1\). Given a word
\(\omega=(\omega_0,\ldots,\omega_{k-1})\in\{1,\ldots,d\}^k\), set \(\Delta_0=[0,1]\).
At step \(r\), the path \(f^{rN}\circ\gamma\) crosses \(R_1\) horizontally
on \(\Delta_r\). By the crossing property of \(F_{\omega_r}\), there is a
nonempty closed subinterval \(\Delta_{r+1}\subset\Delta_r\) on which this path
lies in \(F_{\omega_r}\) and its \(f^N\)-image crosses \(R_1\) horizontally.
After \(k\) steps, choose \(t\in\Delta_k\) and set \(z_\omega=\gamma(t)\).
Since \(t\in\Delta_{r+1}\) for every \(0\le r<k\), we have
\[
 f^{rN}(z_\omega)\in F_{\omega_r}\qquad(0\le r<k).
\]
If two words differ at position \(r\), their points are at least \(\delta\)
apart at time \(rN\). Thus for every \(k\) the \(d^k\) points form a
\((k,\delta/2)\)-separated set for \(f^N\). Since \(\delta\) is independent
of \(k\), taking the exponential growth rate as \(k\to\infty\) gives
\(\log d\le\htop(f^N)=Nh\) by the separated-set definition of entropy.
Thus every finite subfamily has
at most \(e^{Nh}\) members. Substituting \(N=mL+u+v\) gives
\begin{equation}\label{eq:synchronous-upper}
 N_m\le \exp\big((mL+u+v)h\big).
\end{equation}

\subsection{Comparing the entropy bounds}

Let \(J\) be a random variable uniformly distributed on \(J_L\).
Its value \(j\) specifies the split \((j,L-j)\) and hence the product map
\(G_j=f^j\times f^{L-j}\). Independently of \(J\), choose \((x,y)\in D\)
according to \(\rho_D\) and set \(\tau=\tau_J(x,y)\). Write \((\tau,S,T)\)
for the label of the \(\mathcal P_J\)-atom containing \((x,y)\).
The uniform choice gives \(H(J)=\log|J_L|\).
For \(X=(J,\tau,S,T)\) and \(\ell=(j,m,V,W)\), the Shannon entropy is
\[
 H(X)=-\sum_{\ell:\,p_\ell>0}p_\ell\log p_\ell,\qquad
 p_\ell=\Pr(X=\ell)=\frac{\rho_D\bigl(\{\tau_j=m\}\cap(V\times W)\bigr)}{|J_L|}.
\]
By independence, given \(J=j\), each label \((j,m,V,W)\) has conditional
probability equal to the \(\rho_D\)-measure of its corresponding
\(\mathcal P_j\)-atom. Hence \(H(X\mid J=j)=H_{\rho_D}(\mathcal P_j)\).

Conditioning on \(J\) gives, for each \(m\ge1\),
\[
 p_m=\Pr(\tau=m)
 =\frac{1}{|J_L|}\sum_{j\in J_L}\rho_D\bigl(\{\tau_j=m\}\bigr).
\]
Thus the mean return time satisfies, by \eqref{eq:synchronous-entropy},
\[
 \bar\tau=\sum_{m\ge1}mp_m
 =\frac1{|J_L|}\sum_{j\in J_L}\int_D\tau_j\,d\rho_D\le 1/\rho(D).
\]
Set \(q_m=2^{-m}\) for \(m\ge1\). Gibbs' inequality for the probability
distributions \((p_m)\) and \((q_m)\) gives
\[
 H(\tau)\le -\sum_{m\ge1}p_m\log q_m
 =\bar\tau\log 2\le\frac{\log 2}{\rho(D)}.
\]

Given \(\tau=m\) with \(p_m>0\), \(X\) has at most \(N_m\) possible values
across all splits, so
\(H(X\mid\tau=m)\le\log N_m\le(mL+u+v)h\) by \eqref{eq:synchronous-upper}.

Averaging \eqref{eq:label-lower} over \(j\) gives \(H(X\mid J)\ge Lh\bar\tau\).
Since \(J\) and \(\tau\) are coordinates of \(X\), we can apply the entropy
chain rule by conditioning on each in turn. Using the bounds above gives
\[
 \begin{aligned}
 H(X)&=\log|J_L|+\frac{1}{|J_L|}\sum_{j\in J_L} H(X\mid J=j)
       \ge \log|J_L|+Lh\bar\tau,\\
 H(X)&=H(\tau)+\sum_{m:\,p_m>0}p_m H(X\mid\tau=m)
       \le H(\tau)+Lh\bar\tau+(u+v)h.
 \end{aligned}
\]

Comparing the two bounds and cancelling \(Lh\bar\tau\) gives
\[
 \log(L/Q-1)=\log|J_L|\le (u+v)h+H(\tau)
 \le (u+v)h+\frac{\log 2}{\rho(D)} .
\]
The left-hand side tends to infinity with \(L\), while the right-hand side is fixed.
Thus there is at most one ergodic measure of maximal entropy.
By ergodic decomposition, almost every ergodic component of a maximizing
measure also has entropy \(h\). This proves Theorem~\ref{thm:main}.

\subsection*{Acknowledgements}
W. Lin is supported by the National Natural Science Foundation of China
(No.~124B2010).  X. Tian is supported by the National Natural Science
Foundation of China (No.~12471182).

\end{document}